# From the Cherlin-Zilber Conjecture via sharply 2-transitive groups to the Burnside problem

Katrin Tent*

October 13, 2025

**Abstract.** We review the current state of the Cherlin-Zilber Algebraicity Conjecture on simple groups of finite Morley rank, which states that every such group is the group of $K$-rational points of an algebraic group for some algebraically closed field $K$. We will explain the relevance of sharply 2-transitive groups as a potential source of counterexamples and how the Burnside problem necessarily comes into the picture.

**1 Introduction.** The Cherlin-Zilber Algebraicity Conjecture states that any simple group of finite Morley rank is the group of $K$-rational points of an algebraic group for an algebraically closed field $K$. The conjecture is confirmed for groups of Morley rank at most 3 and the general case is still wide open. As in the case of the classification of the finite simple groups the question arises how to determine that an abstractly given group is indeed an algebraic group. For this, the notion of a BN-pair introduced by Tits in the early 1960s will be useful, see Section 3. BN-pairs provide a geometric interpretation of simple Lie groups via the theory of buildings. For split BN-pairs of (Tits-) rank at least 2 and of finite Morley rank the classification is complete, but the case of BN-pairs of rank 1 remains wide open and will be the focus of these notes.

**2 Groups of finite Morley rank** Morley rank is a model theoretic dimension on definable sets which behaves very much like algebraic dimension in the context of algebraic geometry. More precisely it is a rank on formulas of a given (complete and countable) theory $T$ in a fixed first-order language $L$ (see e.g. [TZ]). Note that we generally identify definable sets with the formulas defining them and vice versa. More precisely, given a model $M$ of $T$ and a formula $\varphi(\bar{x}, \bar{a})$ where $\bar{x}$ is an $n$-tuple of variables and $\bar{a}$ is a tuple of elements from $M$, the definable set given by $\varphi(\bar{x}, \bar{a})$ is the set of elements $\bar{b} \in M^n$ such that the formula $\varphi(\bar{b}, \bar{a})$ holds in $M$, i.e.

$$\varphi(M^n, \bar{a}) = \{\bar{b} \in M^n \colon M \models \varphi(\bar{b}, \bar{a})\}.$$

We may also denote this set by $\varphi(M, \bar{a})$ or even $\varphi(M)$ and we often identify the formula $\varphi$ with the set $\varphi(M)$ it defines.

DEFINITION 2.1. *Let $T$ be a complete countable theory with model $M$ and let $\varphi(\bar{x}, \bar{a})$ be a formula with parameters from $M$[1]. Then the Morley rank $MR(\varphi)$ of $\varphi$ is defined inductively as follows:*

1. *$MR(\varphi) \geq 0$ if $\varphi(M) \neq \emptyset$;*
2. *$MR(\varphi) \geq \alpha + 1$ for an ordinal $\alpha$ if there exist formulas $\varphi_i, i < \omega$, such that $\varphi_i(M) \subseteq \varphi(M), i < \omega$, the (sets defined by the) $\varphi_i$ are pairwise disjoint and $MR(\varphi_i) \geq \alpha$;*
3. *$MR(\varphi) \geq \lambda$ for a limit ordinal $\lambda$ if $MR(\varphi) \geq \alpha$ for all $\alpha < \lambda$.*
4. *We put $MR(\varphi) = \alpha$ if the Morley rank is at least $\alpha$, but not at least $\alpha + 1$.*

*If $MR(\varphi) = \alpha$, the* Morley degree *of $\varphi$ is defined as the maximal $k$ such that there exist formulas $\varphi_i, 0 < i < k$ such that $\varphi_i(M) \subseteq \varphi(M), i < k$, the $\varphi_i$ are pairwise disjoint and $MR(\varphi_i) = \alpha$.*

---

*The author was supported by the German Research Foundation (DFG) under Germany's Excellence Strategy EXC 2044–390685587, Mathematics Münster: Dynamics – Geometry – Structure and by CRC 1442 Geometry: Deformations and Rigidity.

[1]We may assume that $M$ is sufficiently saturated, see [TZ].

*We say that the formula $\varphi$ has Morley rank if $MR(\varphi) \not\geq \alpha$ for some ordinal $\alpha$. Otherwise the formula does not have Morley rank.*

*We say that a theory $T$ (or a model $M$ of $T$) has (finite) Morley rank if every formula has (finite) Morley rank or, equivalently, if the formula $\varphi(x) = \ulcorner x = x \urcorner$ has (finite) Morley rank.*

The theories in which every formula has (not necessarily finite) Morley rank are exactly the $\omega$-stable theories (see e.g. [TZ]), a class of theories which is very well-behaved. The class of theories in which every formula has *finite* Morley rank is a much more restrictive class.

Since the Morley rank and Morley degree of a formula can decrease only finitely many times, a group with Morley rank (i.e. any $\omega$-stable group) automatically satisfies the descending chain condition (DCC) on definable subgroups. This implies in particular that any such group has a connected component, namely the smallest definable subgroup of finite index. It is not hard to see that an $\omega$-stable group $G$ is connected if and only if it has Morley degree 1.

The similarity between Morley rank and algebraic dimension can be made precise in the following sense:

Theorem 2.2. *If $X$ is a definable set in an algebraically closed field, then, by the results of Chevalley and Tarski, $X$ is a constructible set and the Morley rank agrees with the Zariski dimension.*

There are a number of results that reinforce the similarity of structures of finite Morley rank to the realm of algebraic geometry, starting with fields:

Theorem 2.3. [Ma] *Any field with Morley rank is algebraically closed.*

On the group side we have the following results:

Theorem 2.4.

1. [Re] *Any group of Morley rank 1 is virtually abelian.*
2. [Ch] *Any group of Morley rank 2 is virtually solvable.*
3. [Fr] *Any group of Morley rank 3 is either virtually solvable or virtually isomorphic to $PSL_2(K)$ for an algebraically closed field $K$.*
4. [ABC] *Groups of finite Morley rank and characteristic* 2 *(in a technical sense) satisfy the Cherlin-Zilber-Conjecture.*

Note that parts 1. − 3. of the previous theorem reflect the same behaviour as algebraic groups of the corresponding dimensions. Parts 1 and 2 can by now be considered classical, whereas part 3. is relatively recent, the first major breakthrough after a long time.

Of course the main question that arises with respect to the Cherlin-Zilber conjecture is how to identify an abstractly given simple group as an algebraic group. In the setting of groups of Morley rank 3, it had been known that unless the group is virtually isomorphic to $PSL_2(K)$ all proper non-trivial definable subgroups are connected and have Morley rank 1. Frécon was able to use these subgroups in order to construct a geometry, which turns out to be (generically) a projective plane, and derive a contradiction. These generic projective planes can be seen as a generalization of the Tits buildings arising from the theory of BN-pairs.

**3 Groups, BN-pairs and associated buildings** Let us start with the purely group theoretic definition of a BN-pair due to Jacques Tits:

Definition 3.1. *Let $G$ be a group, $B, N \leq G$ subgroups. Then $B$ and $N$ form a BN-pair of (Tits)-rank $k$ if the following holds:*

1. $BNB = G$;
2. $T = B \cap N \trianglelefteq N$;
3. *the* Weyl group $W = N/T$ *is a Coxeter group of rank $k$ generated by involutions $s_1, \ldots, s_k$;*
4. $Bs_iBwB = BwB \cup Bs_iwB$ *for $w \in W, i = 1, \ldots, k$.*

*A BN-pair $(B, N)$ splits if there is a normal nilpotent subgroup $U \trianglelefteq B$ such that $B = U.T$.*

*A BN-pair is* spherical *if $W$ is finite.*

Example 3.2. *The standard example of a group with a BN-pair of rank $k$ is the group $G = PGL_{k+1}(K)$ for a field $K$ where we can choose $B$ to be the standard Borel subgroup of $G$, i.e. the subgroup of upper triangular matrices in $G$, and $N$ the normalizer of the subgroup $T$ of diagonal matrices in $G$. Thus $N$ is the subgroup of $G$ consisting of permutation matrices, i.e. matrices having a single entry 1 in each row and column, and 0 otherwise. Then $T = B \cap N$ is a maximal torus and $W = N/T$ is isomorphic to $S_{k+1}$, so is generated by $k$ transpositions. Clearly, by setting $U \trianglelefteq B$ to be the unipotent radical of $B$, i.e. the subgroup of strictly upper triangular matrices with 1s on the diagonal, we see that this BN-pair splits.*

Tits showed that given a BN-pair $(B, N)$ of rank $k$ in a group $G$ there are exactly $2^k$ subgroups $P$ with $B \leq P \leq G$. These subgroups and their conjugates are the *parabolic* subgroups of $G$ and the simplicial complex consisting of all conjugates of parabolics with the edge relation given by inclusion is the associated building.

Abstractly, a building may be defined as a simplicial complex with the following properties:

Definition 3.3. *A building of rank $k$ is a simplicial complex $\Delta$ of rank $k$ together with a class of subcomplexes called apartments such that*

1. *any two $k$-simplices lie in a common apartment;*
2. *for $1 \leq j \leq k-1$ every $j$-simplex of $\Delta$ is contained in at least three $k$-simplices;*
3. *any $(k-1)$-simplex contained in an apartment $A$ lies in exactly two adjancent $k$-simplices of $A$ and the graph of adjacent $k$-simplices is connected;*
4. *if two simplices both lie in apartments $A$ and $A'$, then there is a simplicial isomorphism of $A$ onto $A'$ fixing the vertices of the two simplices.*

*A $k$-simplex in $\Delta$ is called a chamber.*

Example 3.4. *The standard example to keep in mind is the case of $k$-dimensional projective space over any field $K$. The vertices of the simplicial complex are the subspaces with edges given by inclusion. A chamber, i.e. a maximal simplex, is simply a maximal flag, and the apartments can be described as the simplicial subcomplexes one obtains by fixing a basis of the space and considering all subspaces spanned by subsets of this basis.*

*Note that from this description of projective space as a building we can recover the building associated to $PGL_{k+1}(K)$ described above in terms of parabolic subgroups by considering the stabilizers of these (flags of) subspaces. These are exactly the parabolic subgroups and the group $B$ is the stabilizer of a maximal flag. Conversely, we can recover projective space from the building given in terms of parabolic subgroups by looking at the subspaces stabilized by a parabolic subgroup. This is an inclusion reversing bijection.*

The previous example shows the tight connection between groups and buildings. More specifically we have the following:

Proposition 3.5. [Ti] *A group $G$ has a BN-pair of rank $k$ if and only if there is a building $\Delta$ of rank $k$ such that $G$ acts transitively on the set of all pairs $(A, C)$ where $A$ is an apartment in $\Delta$ and $C$ is a chamber in $A$.*

*Proof.* If $G$ has a BN-pair of rank $k$, the transitivity on ordered apartments is easy to see. Conversely, if $G$ acts on a building of rank $k$ transitively on ordered apartments, we let $T$ be the stabilizer of an ordered apartment, $N$ the subgroup of $G$ leaving the apartment invariant and $B$ the stabilizer of the marked chamber of this apartment. Then $B$ and $N$ form a BN-pair of rank $k$: properties 1. and 2. can easily be checked and properties 3. and 4. follow from the fact that apartments are Coxeter complexes, thus corresponding to the Weyl group of $G$ (see e.g. [Gr]). ∎

In the classification of the finite simple groups the classification of finite split BN-pairs of rank 2 due to Fong and Seitz [FS] was an important ingredient. They showed that a split BN-pair of rank 2 in a finite group is the standard BN-pair in a finite group of Lie-type, i.e. $B$ is the Borel subgroup and $N$ the normalizer of a maximal torus in $N$. Their work relied on the classification of finite BN-pairs of rank 1. However, the results about finite groups with BN-pairs of rank 1 have recently turned out to not extend to the infinite setting, see Section 4. Whether corresponding results hold under finite Morley rank assumptions is still open, see Section 2

Despite the lack of corresponding results in the rank 1 case, using a more geometric approach, the results of Fong and Seitz were extended in [TVM, Te04] to infinite groups. To specialize to the rank 2 case, we note that for spherical buildings of rank 2 the definition of a building reduces to that of a generalized $n$-gon:

Definition 3.6. *A generalized $n$-gon is a bipartite graph $\Gamma$ of diameter $n$, girth $2n$ and such that every vertex has valency at least 3. An apartment in $\Gamma$ is a $2n$-cycle and a chamber is an edge.*

Note that for $n = 3$ a generalized 3-gon is exactly a projective plane where the bipartition arises by considering the vertices as points and lines, respectively.

In the finite case, Fong and Seitz could use a fundamental result by Feit and Higman which says that the class of finite generalized $n$-gons is very restricted:

Theorem 3.7. [FH] *If $\Gamma$ is a finite generalized $n$-gon, then $n \in \{3, 4, 6, 8\}$.*

In light of Proposition 3.5 and the Cherlin-Zilber Conjecture one might have hoped that a result similar to Theorem 3.7 might hold for infinite generalized $n$-gons, at least under the assumption of finite Morley rank. However, we have

Theorem 3.8. [Te00b] *For all $n \geq 3$ there exist almost strongly minimal $n$-gons whose automorphism groups act transitively on ordered $2n$-cycles, hence contain BN-pairs of rank 2.*

While it follows from the classification of finite simple groups that in finite groups all BN-pairs split, the corresponding result does not hold for infinite groups, not even if the groups are simple and the underlying building has finite Morley rank:

Theorem 3.9. [GT] *The automorphism groups of the examples constructed in Theorem 3.8 are simple groups with a non-split BN-pair of rank 2.*

However, these automorphism groups are not definable inside the geometries and do not have finite Morley rank. This shows that without further assumptions the general results for infinite groups will be much weaker. However, using the classification of Moufang buildings in [TW] one can extend the classification of split BN-pairs from finite groups to arbitrary groups:

Theorem 3.10. [TVM, Te04, DHTVM] *The following are equivalent for any group $G$:*

1. *$G$ has a split BN-pair of rank $k \geq 2$;*
2. *$G$ is isomorphic to the group $\mathcal{G}(K)$ of $K$-rational points of some algebraic group scheme $\mathcal{G}$ of relative $K$-rank $k$ for some (not necessarily algebraically closed) field $K$.*

*Furthermore, in this case the BN-pair is the standard one given by the Borel subgroup $B$ and the normalizer $N$ of a maximal torus $T$ of $G$.*

Thus, the purely group theoretic notion of a split BN-pair of rank $k \geq 2$ captures the essence of an algebraic group of $K$-rank $k \geq 2$.

As a corollary of the previous theorem we obtain:

Corollary 3.11. *If $G$ is a simple group of finite Morley rank with a split BN-pair of rank at least 2, then $G$ is isomorphic to an algebraic group over an algebraically closed field.*

*Proof.* By Theorem 3.10 the group is algebraic over some field $K$. By [KRT] or more generally [ST] the field $K$ is definable in the abstract group $G$ and hence is algebraically closed by Theorem 2.3. Alternatively, one can deduce from Theorem 3.10 that the associated building satisfies the Moufang condition and thus the result follows from [KTVM]. ∎

We now turn to groups with BN-pairs of rank 1 in the disguise of groups with a 2-transitive group actions. Note that a building of rank 1 is just a set with no further structure.

**4 Sharply 2-transitive groups** Recall that a group action $G$ on a set $X$ is said to be (sharply) $k$-transitive, $k \geq 1$, if for any two $k$-tuples $(x_1, \ldots, x_k), (y_1, \ldots, y_k) \in X^k$ where $x_i \neq x_j, y_i \neq y_j, 1 \leq i \neq j \leq k$ there exists some (unique) element $g \in G$ such that $g(x_i) = y_i, i = 1, \ldots, k$.

*Remark* 4.1.
1. A transitive group action of $G$ on $X$ is (sharply) $n$-transitive if and only if for any $x \in X$, the *stabilizer* $G_x = \{g \in G \colon g(x) = x\}$ acts (sharply) $(n-1)$-transitively on $X \setminus \{x\}$.
2. If the action of $G$ on $X$ is sharply $n$-transitive, then for any distinct elements $x_1, \ldots, x_n \in X$ the stabilizer $G_{x_1,\ldots,x_n}$ is trivial.

3. We call a group $G$ sharply 2-transitive if there exists a set $X$ and an action for $G$ on $X$ which is sharply 2-transitive.

A sharply 1-transitive action is also called *regular* and any regular action arises from right multiplication by a group on itself – clearly there are no restrictions on the groups. The situation changes drastically when one looks at higher degrees of transitivity, see e.g. [Te16] for more background.

The following is well-known:

LEMMA 4.2. *A group $G$ has a BN-pair of rank 1 if and only if $G$ acts 2-transitively on some set $X, |X| \geq 2$.*

*Proof.* If $G$ acts 2-transitively on a set $X$, let $x, y \in X, x \neq y$, let $B = G_x$ be the point-stabilizer of $x$ and $N = G_{x,y}$ the set-wise stabilizer of the set $\{x, y\}$. Then $T = B \cap N = G_{x,y}$, and $N/T \cong S_2$. The other axioms are easily checked.

Conversely, if $G$ has a BN-pair of rank 1, then $B$ is a maximal subgroup in $G$ and $W = N/T \cong S_2$. From $G = BNB$ we see that $G$ acts 2-transitively on the cosets of $B$. □

Here are the standard examples of sharply 2- and 3-transitive groups:

**The case** $n = 2$. Let $K$ be a (not necessarily commutative) field and consider the group $G = AGL(1, K) \cong K \rtimes K^*$ of affine linear transformations $x \mapsto ax + b, a, b \in K, a \neq 0$, which acts sharply 2-transitively on (the affine line of) $K$.

**The case** $n = 3$. Let $K$ be a commutative field. The group

$$G = PGL_2(K) = \{x \mapsto \frac{ax + b}{cx + d} \colon a, b, c, d \in K, ac - bd \neq 0\}$$

acts sharply 3-transitively on the projective line of $K$.

Note that the stabilizer of $\infty$ in $PGL_2(K)$ is isomorphic to $AGL_1(K)$.

**The case** $n \geq 4$. Here the restrictions are even more severe: Jordan proved in 1872 [Jo] that apart from the symmetric and alternating group of degree $n, n + 1$ and $n + 2$, respectively, the only finite sharply $n$-transitive groups are the Mathieu groups $M_{11}$ and $M_{12}$, which are sharply $n$-transitive for $n = 4$ and 5, respectively.

This was generalized by J. Tits and M. Hall, who proved that there are no infinite sharply $n$-transitive groups for $n \geq 4$.

Finite sharply 2- and 3-transitive groups were classified by Zassenhaus in [Z1] and [Z2] in the 1930's and were shown to arise from so-called nearfields, see below. They essentially look like the groups of affine linear transformations $x \mapsto ax + b$ or Moebius transformations $x \mapsto \frac{ax+b}{cx+d}$, respectively.

However, it remained an open problem whether a classification similar to the one in the finite situation holds for infinite sharply 2- and 3-transitive groups.

Note that $G = AGL(1, K)$ has a regular normal subgroup $A \cong (K, +)$, i.e. a normal subgroup acting regularly on the underlying set, consisting of the translations $x \mapsto x + b, b \in K$, and we see that the point stabilizer $G_0$ of $0 \in K$ is the group of homotheties $x \mapsto ax, a \neq 0$, which is obviously isomorphic to the multiplicative group $(K^*, \cdot)$ of the underlying field.

Since $A$ is normal in $G$, the elements of the point stabilizer $G_0$ act by conjugation on the elements of $A$. The action of any $g \in G_0$ on $A$ is given by an automorphism of $A$. In other words, writing the group $A$ additively, for all $n_1, n_2 \in A, g \in G_0$ we have

$$(n_1 + n_2)^g = n_1^g + n_2^g.$$

Thus, even without knowing that the group $G = AGL(1, K) \cong A \rtimes G_0$ arises from a field, the right distributivity appears automatically as a consequence of the group action. For this distributivity, we do not even have to assume that the regular normal subgroup in $G$ is abelian: if a sharply 2-transitive group $G$ (acting on a set $X$) contains a regular normal subgroup $N$, then for any $x \in X$ the point stabilizer $G_x$ acts by conjugation on the normal subgroup $N$ and since the sharp 2-transitivity implies that $N \cap G_x = 1$, we necessarily have

$$G \cong N \rtimes G_x.$$

Now for $g \in G_x, n_1, n_2 \in N$ the action of $g$ on $N$ satisfies

$$(n_1 \cdot n_2)^g = n_1^g \cdot n_2^g.$$

The point stabilizer $G_x$ acts by conjugation as a group of automorphisms of $N$, which is regular on the nontrivial elements of $N$. This already implies (among other things) that all nontrivial elements of $N$ have the same order.

Thus, any sharply 2-transitive group $G$ with a regular normal subgroup gives rise to a *nearfield* in the sense of the following definition:

DEFINITION 4.3. *A (right) nearfield is a structure $Q$, together with two binary operations, $+$ (addition) and $\cdot$ (multiplication), satisfying the following axioms:*

1. $(Q, +)$ *is a group with identity element* $0$ *(we do not need to assume the group to be abelian);*
2. $(Q \setminus \{0\}, \cdot)$ *is a group with identity element* $1$*;*
3. $(a + b) \cdot c = a \cdot c + b \cdot c$ *for all elements* $a, b, c \in Q, c \neq 0$ *(The right distributive law).*

It is easy to see that any nearfield $Q$ gives rise to a sharply 2-transitive group exactly in the same way in which $AGL(1, K)$ arises from a field $K$ and that this group will again have a regular (and in fact abelian) normal subgroup isomorphic to the additive group of $Q$. Namely for $a, b \in Q, a \neq 0$, the group of transformations

$$\varphi_{a,b} \colon Q \longrightarrow Q, x \mapsto a \cdot x + b$$

acts sharply 2-transitively on $Q$. In this way, there is a one-to-one correspondence between nearfields and sharply 2-transitive groups having a regular normal subgroup.

*Remark* 4.4. It is easy to see that in any 2-transitive group every nontrivial normal subgroup must be transitive and that transitive abelian groups act regularly. Thus, a 2-transitive group $G$ with a non-trivial abelian normal subgroup $A$ splits as $G = A \rtimes G_x$.

DEFINITION 4.5. *A sharply 2-transitive group $G$ splits if it contains a nontrivial abelian normal subgroup $A$, so $G \cong A \rtimes G_x$.*
*Equivalently, $G$ splits if it arises from a nearfield $K$ as $G \cong K \rtimes K^*$.*
*We call a sharply 2-transitive group $G$* standard *if $G \cong K \rtimes K^*$ for some algebraically closed field $K$.*

Clearly, the existence of a nontrivial abelian normal subgroup is independent of the particular action of the group on the underlying set.

In the finite case, Zassenhaus proved that every sharply 2-transitive group arises from a nearfield and classified all finite nearfields by characterizing those finite (linear) groups that can act regularly on the set of non-zero vectors of a finite vector space. It is interesting to note that while there are no finite (proper) skew fields, proper finite nearfields do exist.

*Remark* 4.6. **A word of caution:** In general, the two notions of splitting introduced above in the context of sharply 2-transitive groups do not agree:

(i) The BN-pair arising from a sharply 2-transitive group splits in the sense of Definition 3.1 if and only if the point-stabilizer $B$ is nilpotent (note that the stabilizer $T$ of two distinct points is trivial).

(ii) A sharply 2-transitive group $G$ is split in the sense of Definition 4.5 if the point-stabilizer $B = G_x$ has an abelian complement $A$.

Note that in a split sharply 2-transitive group, point stabilizers need not be nilpotent as can be seen by considering $AGL_1(k)$ for a skewfield $k$, hence (ii) does not imply (i).

It is still unknown whether (i) implies (ii). So far there are no non-split sharply 2-transitive groups with nilpotent point-stabilizers. Note however that under finite Morley rank assumptions, these two notions converge, see Corollary 5.3.

In many classes of groups, classification results analogous to the finite case hold. Before stating these results we introduce the characteristic of a sharply 2-transitive group:

**4.1 The characteristic of a sharply 2-transitive group** It is easy to see that any sharply 2-transitive group $G$ acting on a set $X$ contains plenty of involutions and that the set

$$J = \{t \in G \setminus \{1\} : t^2 = 1\}$$

of involutions in $G$ forms a single conjugacy class. Hence either all involutions of $G$ have a fixed point (and if they do each involution has a unique fixed point) or no involution has a fixed point.

The following important observation is well-known (for a proof see e.g. [Te16]).

LEMMA 4.7. *If $G$ is sharply 2-transitive on $X$ and involutions in $G$ have fixed points, there is a natural $G$-equivariant bijection*

$$J \longrightarrow X, t \mapsto Fix(t)$$

*where $Fix(t)$ denotes the fixed point of $t$ and the action of $G$ on $J$ is given by conjugation.*

Thus in the case where involutions have fixed points we have the following consequences:

COROLLARY 4.8. *If $G$ is sharply 2-transitive on $X$ and involutions in $G$ have fixed points, then*

1. *every point stabilizer $G_x$ contains a unique involution (and so this involution is central in $G_x$);*
2. *the set of products of two distinct involutions forms a conjugacy class.*

The second part follows since an element $g \in G$ taking the pair $(x_1, y_1)$ of distinct elements of $X$ to the pair $(x_2, y_2)$ will conjugate the (ordered) pair of involutions with fixed points $x_1, y_1$ to the pair with fixed points $x_2, y_2$.

If involutions have no fixed points, then not only does the previous argument not work, but in fact we will see in Theorem 4.13 below that $J^2 \setminus \{1\}$ does not necessarily form a single conjugacy class.

DEFINITION 4.9. *A* translation *in a sharply 2-transitive group is a product of two distinct involutions, i.e. an element of $J^2 \setminus \{1\}$.*

Note that the set $J^2$ is not necessarily a subgroup!

In the group $AGL(1, K)$ the involutions are exactly the elements of the form $x \mapsto -x + b$. Thus in this case the set

$$J^2 = \{t_1 t_2 : t_1, t_2 \in J\}$$

equals the set of translations $x \mapsto x + b, b \in K$, which is isomorphic to the additive group of the field $K$. Note that if the field $K$ is of characteristic 2, then the involutions $x \mapsto -x + b$ are fixed point free (and commute pairwise). This explains and motivates the following definition

DEFINITION 4.10. *Let $G$ be a sharply 2-transitive group.*

1. *If involutions in $G$ have fixed points, we say that $G$ has characteristic $p$ if the elements of $J^2 \setminus \{1\}$ have order $p$ and characteristic 0 if their order is infinite.*
2. *If involutions in $G$ are fixed point free, we say that $G$ has characteristic 2.*

Note that the characteristic of a sharply 2-transitive group is necessarily a prime $p$ and the case $p = 2$ occurs if and only if involutions are fixed point free (as otherwise two distinct involutions cannot commute). Note also that if the characteristic of $G$ is 2, then the involutions in $G$ commute if and only if $J^2$ is a (commutative) subgroup.

A result of B. Neumann [Ne] now gives a partial answer to the above question or at least a clear criterion:

THEOREM 4.11. *If $G$ acts sharply 2-transitively on the set $X$, then $G$ contains a regular normal subgroup $N$ if and only if for some (or any) involution $t \in G$ the set*

$$tJ := \{tt_1 : t_1 \in J\}$$

*is a subgroup. In this case we have $N = tJ = J^2$ and the group is abelian.*

Note that Neumann's result says in particular that addition in nearfields is always commutative. Thus, any regular normal subgroup of a sharply 2-transitive group is abelian. Hence Definition 4.5 is equivalent to saying that a sharply 2-transitive group splits if it contains a *regular* normal subgroup.

Until recently it was widely believed (and conjectured) that all sharply 2-transitive groups split. Here is a list of splitting results in different settings:

Theorem 4.12. *Let $G$ be a sharply 2-transitive group. Then $G$ splits in any of the following cases:*

1. [Z1] *$G$ is finite;*
2. [Ti52, Ti56] *$G$ is locally compact connected;*
3. [W] *$G$ is locally finite;*
4. [Te00a] *$G$ is definable in an o-minimal structure;*
5. [GlGu1] *$G$ has characteristic $\neq 2$ and is linear over a field of characteristic $\neq 2$;*
6. [GMS] $\mathrm{char}(G) \neq 2$ *and if (i)* $\mathrm{char}(G) \neq 0$*, every subgroup generated by three involutions is either linear or solvable and if (ii)* $\mathrm{char}(G) = 0$*, $G$ is linear;*
7. [Ke, Thm. 9.5, p. 42] $\mathrm{char}(G) = 3$ *(see also* [Tu]*).*
8. [M] *The exponent of the point stabilizer is 3 or 6;*

However, by now there is a vast literature on non-split sharply 2-transitive groups. The first such groups were constructed by Rips, Segev and Tent [RST]:

Theorem 4.13. [RST] *Any group embeds into a non-split sharply 2-transitive group in characteristic 2.*

In particular, in such groups dihedral subgroups, i.e. subgroups generated by two distinct involutions, need not be conjugate.

In the case of characteristic different from 2, the action of a sharply 2-transitive group is necessarily the conjugation action on the set of involutions by Lemma 4.7. Thus, in order to obtain a sharply 2-transitive group we can proceed by iterated HNN-extensions and obtain:

Theorem 4.14. [RT] *Any torsion free group embeds into a non-split sharply 2-transitive group in characteristic 0.*

**4.2 Linearity** A different version of Theorem 4.13 was published in [TZ]. The paper shows in particular that the group $(\mathbb{Z}/2\mathbb{Z} \times F_\omega) * F_\omega$ is a non-split sharply 2-transitive group of characteristic 2, where $F_\omega$ is the free group on countably many generators.

As pointed out by Simon André, this group embeds into $(\mathbb{Z}/2\mathbb{Z} \times F_2) * F_2$ where $F_2$ is the free group on two generators. Since $(\mathbb{Z}/2\mathbb{Z} \times F_2) * F_2$ is virtually free, this group is virtually linear, hence linear.

This shows that in Theorem 4.12 5. the assumption that the permutation characteristic is different from 2 is necessary. Note that [GlGu2] also construct sharply 2-transitive groups that embed into $SL_n(\mathbb{R})$ for $n \geq 3$.

**4.3 Sharply 3-transitive groups** It also turns out that at least in characteristic 2, non-split sharply 2-transitive groups can appear as point-stabilizers of sharply 3-transitive groups, thus showing that the class of such groups is much wilder than expected:

Theorem 4.15. [Te2] *There exist non-split sharply 3-transitive groups of characteristic 2.*

Here, the characteristic of a sharply 3-transitive group is defined as the characteristic of the (sharply 2-transitive) stabilizers (see Remark 4.1) and such a group is called non-split if the point stabilizers are non-split sharply 2-transitive groups.

**Question** Are there non-split sharply 3-transitive groups of characteristic different from 2?

There seem to be serious restrictions on such constructions, but very little is known. In contrast to the sharply 2-transitive case and the sharply 3-transitive groups without further assumptions, sharply 3-transitive groups of finite Morley rank have long been classified. Nesin's proof relies on Hrushovski's fundamental work showing that any definable transitive group action on a strongly minimal set is essentially sharply $k$-transitive for $k \in \{1, 2, 3\}$ [Hr]:

Theorem 4.16. [Ne] *Every sharply 3-transitive group of finite Morley rank is isomorphic to the group $PSL_2(K)$ for an algebraically closed field $K$ acting on the projective line of $K$ in the natural way.*

**5 Sharply 2-transitive groups of finite Morley rank** In the setting of finite Morley rank it is still open whether any sharply 2-transitive group splits. Here is the list of known results:

THEOREM 5.1. *Let $G$ be a sharply 2-transitive group of finite Morley rank, then $G$ splits in any of the following cases:*

1. [ABW] $\mathrm{char}(G) = 2$*;*
2. [BN] (Thm. 11.73) $G_x$ *is nilpotent;*
3. [CT] $G$ *has Morley rank* 6*.*

Furthermore, nearfields of finite Morley rank have recently been classified. In contrast to the finite case, no proper nearfields of finite Morley rank exist, at least in characteristic different from 2:

THEOREM 5.2. [ABW] *Every nearfield of finite Morley rank and characteristic different from* 2 *is an algebraically closed field.*

In light of Theorem 5.1 we therefore obtain

COROLLARY 5.3.

1. *Any split sharply* 2*-transitive group of characteristic* $\neq 2$ *and of finite Morley rank is standard.*
2. *For sharply* 2*-transitive groups of finite Morley rank and characteristic* $\neq 2$ *the two notions of splitting (see Remark 4.6) coincide.*
3. *Any sharply* 2*-transitive group of Morley rank* 6 *is standard.*

*Remark* 5.4. It is well-known that if the Cherlin-Zilber Conjecture holds, then every sharply 2-transitive group of finite Morley rank is standard. In order to prove directly that every sharply 2-transitive group of finite Morley rank is indeed standard, in light of Corollary 3.11 it is thus left to prove that

1. any sharply 2-transitive group of finite Morley rank and characteristic different from 2 or 3 splits; and
2. any nearfield of finite Morley rank and characteristic 2 is an algebraically closed field.

On the other hand, certain non-split sharply 2-transitive groups potentially indeed yield counterexamples to the conjecture.

PROPOSITION 5.5. [CT] *A nonsplit sharply* 2*-transitive group of finite Morley rank in which the centralizers of translations are strongly minimal is simple. Hence such a group would be a counterexample to the Cherlin-Zilber-Conjecture.*

Note that by [AT] and [AG] there do exist simple sharply 2-transitive groups!

Let $G$ be a sharply 2-transitive group of characteristic $\mathrm{char}(G) \neq 2$ and let $J$ be the set of involutions. If commuting is transitive on the set of non-trivial translations in $J^2$, there is a well-behaved point-line geometry first defined by Schröder [Sch] following ideas of Bachmann [Ba].

LEMMA 5.6. [BN](Lemma 11.50) [CT] *Let $G$ be a sharply* 2*-transitive group of characteristic* $\neq 2$*. Then the following conditions are equivalent:*

*(a) Commuting is transitive on* $J^2 \setminus \{1\}$*.*

*(b)* $iJ \cap kJ$ *is uniquely 2-divisible for all involutions* $i \neq k \in J$*.*

*(c)* $Cen(ik) = iJ \cap kJ$ *is abelian and is inverted by* $k$ *for all* $i \neq k \in J$*.*

*(d) The set* $\{Cen(\sigma) \setminus \{1\} : \sigma \in J^2 \setminus \{1\}\}$ *forms a partition of* $J^2 \setminus \{1\}$*.*

*If $G$ has finite Morley rank or* $\mathrm{char}(G) \neq 2, 0$*, these conditions are satisfied.*

In the setting of the previous lemma, we can consider elements of $J$ as points and the sets $iJ$ as lines. In fact, this point-line geometry is very helpful in the investigation of sharply 2-transitive groups of finite Morley rank:

Theorem 5.7. [CT] *Let $G$ be a sharply 2-transitive group of finite Morley rank and characteristic $\neq 2$. If the associated point-line geometry is a generic projective plane, then $G$ is standard.*

In [CT] this point-line geometry was extended to a larger class of groups including Frobenius groups using the definition of a mock hyperbolic reflection space:

Definition 5.8. [CT] *Let $G$ be a group and $J \subset G$ a conjugacy class of involutions in $G$. For involutions $i \neq j \in J$ we write $\ell_{ij} = \{k \in J : ij \in kJ\}$ and let $\Lambda = \{\ell_{ij} : i \neq j \in J\}$. Then $(J, \Lambda)$ is a mock-hyperbolic reflection space if the following conditions are satisfied:*

1. *any two points $i, \neq j \in J$ lie on a unique line;*
2. *for any $i, j \in J$ there is a unique $k \in Q$ with $i^k = j$;*
3. *if $\lambda^i = \lambda^j$ for $i \neq j \in J, \lambda \in \Lambda$, then $\lambda^i = \lambda = \lambda^j$.*

Examples of such spaces arise from the hyperbolic plane and $PSL_2(\mathbb{R})$, from sharply 2-transitive groups and more generally Frobenius groups (see [CT]) and from bad groups in the sense of Cherlin. These geometries can be used to detect splitting in a wider context:

Theorem 5.9. [CT] *If $G$ is a group with an associated mock hyperbolic reflection space $J$, then the following are equivalent:*

1. *$G \cong A \rtimes \mathrm{Cen}(i)$ for some abelian normal subgroup $A$ and any $i \in J$;*
2. *$J$ is a (possibly degenerate) projective plane;*
3. *$J$ consists of a single line.*

On the other hand, the examples constructed in [RST] (see also [TZ]) show that in characteristic 2 these conditions need not be satisfied as *any group* can be embedded into a sharply 2-transitive group of characteristic 2. However, the non-split examples in characteristic 0 constructed in [RT] do satisfy the assumptions and it is an open question whether non-split sharply 2-transitive groups exist in characteristic 0 which fail to satisfy these conditions.

Until recently all known sharply 2-transitive groups were split and so there seemed little point in looking for simple groups of finite Morley rank in this context. However, in light of the new examples of sharply 2-transitive groups, it becomes more relevant to investigate simple sharply 2-transitive groups.

*Remark* 5.10. [CT] Note that by Zilber's Indecomposability Theorem, in a sharply 2-transitive group of finite Morley rank there exists a bound $m$ such that every element can be written as a product of $m$ involutions.

This criterion suffices to show that the groups constructed in [RST] and in [RT] do not have finite Morley rank. The groups constructed in [AT] satisfy this condition with $m = 4$, the minimal possible value for $m$. But by an observation of Frank Wagner these groups do not have finite Morley rank:

Proposition 5.11. *Suppose $G$ is a simple sharply 2-transitive group, $t \in G$ is a translation with infinite centralizer and $G_0 = Cen(t) \lneq G_1 \lneq \ldots \lneq \bigcup_{i<\omega} G_i = G$. Then $G$ does not have have finite Morley rank.*

*Proof.* Suppose otherwise and let $C = Cen(t)$ denote the centralizer of $t$ in $G$. Then $C$ is connected by [BN] (Lemma 11.60) and by Zilber's Indecomposability Theorem, the subgroups $N_i = \langle C^g : g \in G_i \rangle \trianglelefteq G_i$ are definable and connected. Thus, this ascending chain of $N_i, i < \omega$, becomes stationary and its maximal element is a normal subgroup in $G$. □

Since the groups constructed in [AT, GT] satisfy the conditions of the previous proposition, they do not have finite Morley rank. The previous proposition shows one reason why it is difficult to construct simple sharply 2-transitive groups of finite Morley rank and characteristic 0 and thus puts the focus on positive characteristic, where the order of a translation is finite.

Note that the examples constructed in [AAT] in positive characteristic show that there exist infinite sharply 2-transitive groups in which the centralizer of a translation is finite, see Theorem 6.5.

**6 Constructing sharply 2-transitive groups of positive characteristic** As explained in Section 4, in the case of a sharply 2-transitive group of characteristic $\neq 2$ the problem is entirely group theoretic: A group $G$ acts sharply 2-transitively on a set $X$ in characteristic different from 2 if and only if the group acts sharply 2-transitively by conjugation on its set of involutions.

Suppose we are given a group $G$ containing some involutions and such that all dihedral subgroups, i.e. subgroups generated by two distinct involutions, are infinite. Then we want to extend the group $G$ to a group in which all pairs of distinct involutions are conjugate. In the first step we need to ensure that the involutions form a single conjugacy class. The natural method to make elements or subgroups of a given group $G$ conjugate is by HNN-extensions. These are extensions of $G$ by a new element, the stable letter $t$ with the new relation enforcing that $t$ conjugates the initial group to the intended target group. Given the group $G$ and involutions $i, j \in G$ the group

$$H = \langle G, t \mid i^t = j^t \rangle$$

is an extension of $G$ in which $i, j$ are conjugate under $t$. Note that in such an extension any element of finite order in $H$ is conjugate to an element of the group $G$. Thus, the new involutions that appear in this extension $H$ are already conjugate to the existing ones in $G$.

Thus we can use iterated HNN-extensions to eventually ensure that all involutions are conjugate, i.e. that the new group $G'$ acts transitively by conjugation on the set of involutions. In the next step, for a sharply 2-transitive action we need to make the point-stabilizer, in this case the centralizer of an involution $i$, transitive on the set of remaining involutions, see Remark 4.4.

Suppose that $j$ and $k$ are involutions in our new group $G'$ not yet conjugate under the centralizer of $i$. Then we use the HNN-extension

$$H' = \langle G', t \mid i^t = i, j^t = k \rangle$$

to add an element $t$ to the centralizer of $i$ conjugating $j$ to $k$.

Note that $H'$ contains new involutions, which are conjugate to involutions in $G'$, but not necessarily under elements centralizing $i$. By an infinite chain we can assure that finally in the group $\widetilde{G}$ obtained as the limit the conjugation action is sharply 2-transitive on involutions.

This procedure works well in characteristic 0 and 2 because the translations appearing in the HNN-extensions will have infinite order and in characteristics 0 and 2 this is fine.

However, for odd positive characteristic we need to take quotients that ensure that the translations have (the correct) finite order and the result of Kerby shows that this is not possible for $p = 3$.

It is natural to consider a quotient of the group $\widetilde{G}$ of the form

$$\widehat{G} = \widetilde{G} / \langle (ij)^p : i, j \in \widetilde{G} \rangle$$

in which all translations have order $p$. However, it is generally a difficult problem to establish that the quotient group is still infinite and not split. In fact, we know that for $p = 3$ this is impossible by Theorem 4.12.7! The problem of controlling the quotient is thus reminiscent of Burnside's question from 1902 whether every finitely generated torsion group is finite.

**6.1 Small cancellation conditions and Burnside groups** While there have been many negative answers to Burnside's question (see the literature cited in [ART]) we now turn to the most recent approach to showing the infinity of such groups for sufficiently large exponent. For a subset $R$ of a group $G$ we write $\langle\langle R \rangle\rangle$ to denote the normal subgroup of $G$ generated by the elements in $R$.

Definition 6.1. *The free Burnside group $B(m, n)$ on $m$ generators and exponent $n$ is the quotient of the free group $F_m$ of rank $m$ by all $n^{th}$ powers, i.e.*

$$B(m, n) = F_m / \langle\langle w^n : w \in F_m \rangle\rangle.$$

For exponents $n = 2, 3, 4$ and 6 it is known by work of Burnside, Sanov and M. Hall, respectively, that the free Burnside group is indeed finite for any finite number $m$ of generators. There exist a number of results due among others to Adian and Novikov, Olshanskii, Gromov, Lysenok and Coulon, showing that for sufficiently large (odd) exponents $n$ the free Burnside groups are indeed infinite. However, these results either need an extremely large exponent or the proofs require an extremely large number of induction hypotheses. The situation for even exponents is even more complicated.

So in order to construct sharply 2-transitive groups of sufficiently large positive characteristic, we needed a more direct proof and a better exponent to be able to adapt it to the sharply 2-transitive setting.

We used iterated small cancellation conditions to prove:

THEOREM 6.2. [ART] *The free Burnside group $B(m,n)$ is infinite for odd $n \geq 557$ and any $m \geq 2$.*

This is the lowest currently known lower bound for the exponents of infinite free Burnside groups. The proof is based on iterated small cancellation methods which we will explain next.

**6.2 Small cancellation conditions** Even though the free Burnside groups are not hyperbolic (every element is torsion!) the construction shares many features with that of hyperbolic groups. Recall that a finitely generated group $G = \langle X \rangle$ is $\delta$-hyperbolic in the sense of Gromov (for some $\delta \geq 0$) if all geodesic triangles in the Cayley graph $\Gamma(G, C)$ are $\delta$-thin, i.e. if every side of the triangle is contained in the $\delta$-neighborhood of the other two sides. It is not hard to see that for finitely generated groups being hyperbolic does not depend on the generating set.

Hyperbolicity has strong implications on the algebraic structure of the group. For example using the fact that triangles are thin, it is not hard to see that a hyperbolic group $G$ generated by a finite set $X$ is *finitely presented* (see Proposition 6.3). This means that there is a finite set of *relators* $R \subseteq F(X)$ in the free group $F(X)$ over $X$ such that

$$G \cong F(X)/\langle\langle R \rangle\rangle =: \langle X \mid R \rangle.$$

Given a group $G$ with a presentation $G = \langle X \mid R \rangle$, a word $w \in F(X)$ in the generators of $G$ represents the identity in $G$ if and only if $w$ can be written as a product of conjugates of relators in $R$, or, equivalently, if the word corresponds to the labels on a cycle in the Cayley graph $\Gamma(G, X)$ of $G$. A word $w$ over the letters of $X \cup X^{-1}$ is called *reduced* if it does not contain a subword of the form $xx^{-1}$ or $x^{-1}x$ for $x \in X$.

Here is an example of how the hyperbolicity of the group influences its algebraic properties:

PROPOSITION 6.3. *If $G = \langle X \rangle$ is a $\delta$-hyperbolic group, $X$ finite, then $G$ has a finite presentation $G = \langle X \mid R \rangle$ where $R$ is the collection of all reduced words over $X \cup X^{-1}$ of length less than $8 \cdot \delta$ which represent the identity in $G$.*

*Proof.* Suppose that $G = \langle X \rangle$ is $\delta$-hyperbolic and let $w \in F(X)$ be a reduced word in the generators of $G$ that represents the identity in $G$ and assume that the length $|w|$ of $w$ is at least $8 \cdot \delta$. Clearly we may assume $\delta \in \mathbb{N}$.

Consider the word $w = x_0 \dots x_m$ as a cycle $\gamma = (x_0, \dots, x_m, x_{m+1} = x_0)$ in the Cayley graph $\Gamma(G, X)$. We claim that there exist $0 \leq i \neq j \leq m$ such that the graph distance $d(x_i, x_j)$ in the Cayley graph is strictly less than the distance between $x_i$ and $x_j$ on the cycle $\gamma$ given by $w$.

Suppose otherwise and write $m = 2k$ or $m = 2k + 1$, depending on the parity of $m$. Then the paths $(x_0, \dots, x_k), (x_k, \dots, x_{k+2\delta})$ and $(x_{k+2\delta}, \dots, x_m)$ are geodesic and hence form the sides of a geodesic triangle $\Delta(x_0, x_k, x_{k+2\delta})$. By $\delta$-hyperbolicity, the geodesic $(x_0, \dots, x_k)$ is in the $\delta$-neighbourhood of $(x_k, \dots, x_{k+3\delta}) \cup (x_{k+3\delta}, \dots, x_m)$. Hence the vertex $x_{2\delta} \in (x_0, \dots, x_k)$ is in the $\delta$-neighbourhood of the path $(x_k, \dots, x_m)$. Thus, there exists some $k \leq i \leq m$ such that $d(x_{2\delta}, x_i) \leq \delta$, contradicting our assumption.

This shows that any word $w \equiv 1_G$ of length $|w| > 8\delta$ can be written as a product of words of shorter lengths, each representing the identity in $G$ (see the picture below).

In the picture, assume $w = u_1 u_2 \equiv 1_G$. Then $r_1 = u_1 v, , r_2 = v^{-1} u_2 \equiv 1_G$ and $w = r_1 r_2$:

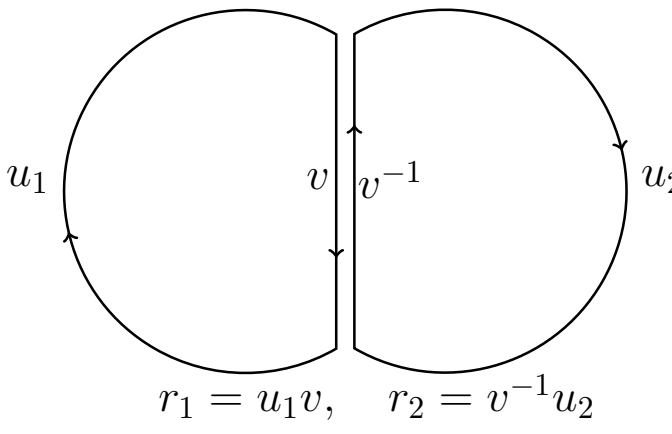

This finishes the proof of the proposition. ∎

The previous proof shows furthermore that every word $w$ that represents the identity in $G$ contains (up to a cyclic shift) a subword $u$ which is contained in one of the relators $r$ and the length of $u$ is strictly larger than half of the length of $r$. In other words we see that the word $w$ contains more than one half of a relator from $R$ as a subword. Presentations with this property are called *Dehn presentations*. Clearly, in groups with a Dehn presentation the word problem is decidable: reduce the length of the word by multiplying with (the inverse of) a relator until the word does not contain any subwords of relators that are strictly longer than half the relator. If in this way the word is reduced to the identity, it presents the identity, otherwise it does not.

THEOREM 6.4. [DK] *A finitely generated group is hyperbolic if and only if it has a Dehn presentation.*

These nice algorithmic properties of hyperbolic groups can – to some extent – also be applied to groups that are not finitely presented as long as the overlaps between different relators (i.e. the pieces of relators that cancel with each other) are well-controlled:

The proof of Theorem 6.2 works inductively by choosing a canonical representative for every coset in $B(m,n) = F_m/\langle\langle w^n : w \in F_m\rangle\rangle$. The induction is based on the rank of a word $w \in F_m$ where (roughly speaking) we define the rank $rk(w)$ to be greater or equal to $k+1$ with respect to our *nesting constant* $\tau = 15$ if the word $w$ (cyclically) contains a subword of the form $v^\tau$ for some word $v \in F_m$ with $rk(v) \geq k$.

We define

$$N_k = \langle w^n : rk(w) \leq k\rangle.$$

Thus, we obtain an ascending sequence of normal subgroups

$$N_0 \leq N_1 \ldots \leq N_i \leq \ldots \bigcup N_i = N = \langle w^n : w \in F_m\rangle.$$

We inductively define the canonical form $can_k(w)$ for a word $w$ as a canonical representative for $wN_k$. In particular, for all $w_0, w_1 \in F_m$ we have

$$w_0 N_k = w_1 N_k \text{ if and only if } can_k(w_0) = can_k(w_1)$$

and we can define a group operation on the set of canonical forms of rank $k$ making this group isomorphic to $F_m/N_k$.

At stage $k$, suppose we are given a word $w$ containing subwords of relators $r_1, \ldots, r_4 \in N_k$.

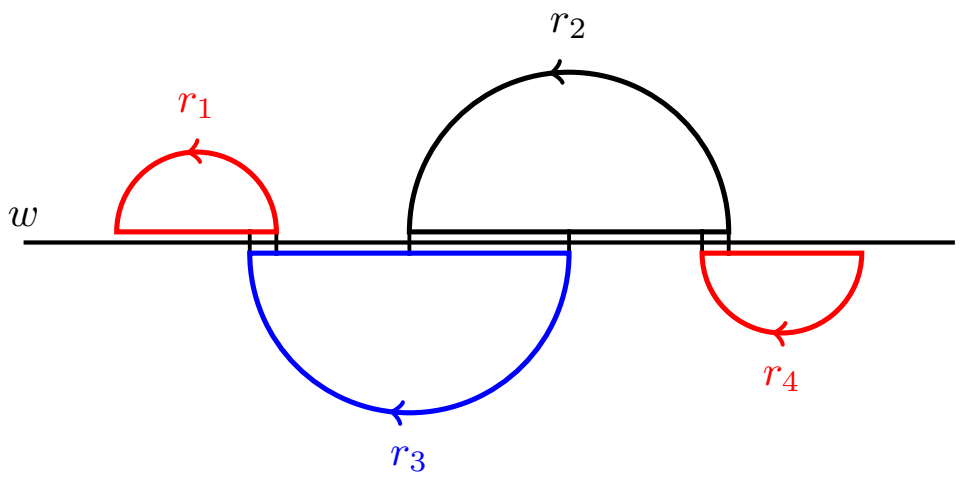

Then we can modify $w$ by choosing any side of these relators without changing the coset of $w$ modulo $N_k$. The small cancellation condition ensures that the overlaps between different relators are relatively short. So e.g. choosing a side for the relator $r_2$ will change the subwords of $r_1$ and $r_3$ only in a relatively small way. To make sure that a choice e.g. for $r_1$ does not lead to problems with the choices for $r_2, \ldots, r_4$ we use the concept of a *certification sequence* that allows us to make an informed decision about a possible turn of a relator that will not lead to problems further away from this particular occurrence.

The important point is that for any $w \in F_m$, the canonical form stabilizes, i.e. for any $w \in F_m$ there is some $k$ such that $can_k(w) = can_l(w)$ for all $l > k$ and thus $can_k(w)$ will be the canonical representative for $wN \in B(m,n)$.

In this way we obtain a section $can : F_m/N \longrightarrow F_m$ i.e. we have

$$can(w) = can(w') \text{ if and only if } wN = w'N \in B(m,n).$$

The set of canonical forms $can(F_m)$ with the appropriate multiplication then forms a group isomorphic to $B(m,n)$.

Thus, the main thrust of [ART] lies in inductively defining $can_k(w)$ for any $k$ based on 13 induction hypotheses. It turns out that any cube-free element of $F_m$ is already in canonical form and so the infinity of the Burnside group follows immediately from the fact that there are infinitely many cube-free words on two letters.

For the methods in [ART] to give a relatively short and accessible proof, we currently need the exponent $n$ to be at least $n > 36 \cdot 15 + 16 = 556$. The proof also yields (the previously known result) that the infinite free Burnside groups are not finitely presented.

**6.3 Sharply 2-transitive groups of positive characteristic** We now turn back to sharply 2-transitive groups. In light of the fact that every sharply 2-transitive group in characteristic 3 splits, it is remarkable that characteristic 3 is not indicative for other characteristics:

THEOREM 6.5. [AAT] *There exist non-split sharply 2-transitive group of prime characteristic $p$ for every sufficiently large $p$.*

In fact, Amelio proves the following stronger statement:

THEOREM 6.6. *[Am] For every sufficiently large prime $p$ congruent 3 mod 4 and all sufficiently large odd numbers $q_1, q_2$ there exists a countable non-split sharply 2-transitive group $G$ of characteristic $p$ and exponent $lcm(q_1, q_2, p, p-1)$. In addition, any $g \in G$ satisfies one of the following:*

- *$g$ is contained in a subgroup of $G$ isomorphic to $AGL(1, \mathbb{F}_p)$;*
- *$g$ is contained in a subgroup of $G$ isomorphic to $\mathbb{Z}/\mathbb{Z}_{q_1}$ if $g$ centralizes no involution;*
- *$g$ is contained in a subgroup of $G$ isomorphic to $\mathbb{Z}/\mathbb{Z}_{q_2}$ if $g$ centralizes an involution.*

*Furthermore, there exist elements in the centralizer of an involution that do have order $q_2$.*

This gives the first examples of sharply 2-transitive in positive characteristic in which point stabilizers have finite exponent, negatively answering the question raised in [M] whether such groups are necessarily split (let alone finite).

The proofs of the previous two theorems rely on Rémi Coulon's results on periodic quotients [Co]. In fact, Coulon uses geometric small cancellation theory, i.e. an analog of the combinatorial small cancellation theory explained above in the context of groups acting on trees. Since the HNN-extensions naturally yield actions on the corresponding Bass-Serre trees, Coulon's results provide the perfect approach to taking the required quotients for the constructions.

*Remark* 6.7. Note that the examples presented in this section do not have finite Morley rank by Remark 5.10. However, as pointed out above, the argument provided in Proposition 5.11 does not apply to these groups as the centralizers of translations are finite.

The examples constructed here also provide a negative answer to Question 11.52 of [BN]: in our examples, the centralizers of translations are cyclic (and of order $p$) and the normalizers are isomorphic to $AGL(1, \mathbb{F}_p)$, hence finite.

**7 Outlook** The last few years have seen huge progress on the structure of sharply 2-transitive groups (with and without finite Morley rank). In the context of the Cherlin-Zilber Conjecture, the urgent questions that remain are:

1. Are there non-split sharply 2-transitive groups of finite Morley rank?

   Note that the existence of such groups would refute the Cherlin-Zilber Conjecture even if the groups are not simple (see Remark 5.4).

2. Is there a simple sharply 2-transitive group of positive characteristic in which every element is a product of a bounded number of translations?

   Such a group would be a candidate for a group of finite Morley rank. Note that in such a group the lines in the associated mock hyperbolic reflection space are finite.

3. Is there a mock hyperbolic reflection space of finite Morley rank in which lines are strongly minimal and in the prime model, every point is definable over any two other points.

For almost one hundred years the class of sharply 2-transitive groups was expected to be extremely rigid. It turned out that in contrast to these expectations, this class is absolutely wild and by now there are non-split sharply 2-transitive groups with many additional properties. What we are still missing is an answer to the question whether there are such groups of finite Morley rank.

**Acknowledgments.** I would like to thank Marco Amelio, Simon André, Anna Cascioli, and Theo Grundhöfer for their comments on earlier versions of these notes and Hannah Kramer for providing the graphics.

Katrin Tent
Mathematisches Institut und
Institut für mathematische Logik
Universität Münster
Einsteinstr. 62
48149 Münster
tent@wwu.de